\documentstyle[11pt,a4wide]{article}
\font\tenmsb=msbm10 scaled \magstep 1
\font\sevenmsb=msbm7 scaled \magstep 1
\font\fivemsb=msbm5 scaled \magstep 1
\newfam\msbfam
\textfont\msbfam=\tenmsb
\scriptfont\msbfam=\sevenmsb
\scriptscriptfont\msbfam=\fivemsb
\def\Bbb#1{\fam\msbfam\relax#1}

\font\teneufm=eufm10 scaled \magstep 1
\font\seveneufm=eufm7 scaled \magstep 1
\font\fiveeufm=eufm5 scaled \magstep 1
\newfam\eufmfam
\textfont\eufmfam=\teneufm
\scriptfont\eufmfam=\seveneufm
\scriptscriptfont\eufmfam=\fiveeufm
\def\frak#1{{\fam\eufmfam\relax#1}}

\title{\bf ON q-ANALOGUES OF BOUNDED SYMMETRIC DOMAINS AND DOLBEAULT
COMPLEXES}

\author{\sl S. Sinel'shchikov \thanks{Partially supported by ISF grant
U2B200 and grant DKNT-1.4/12}\and \hspace{-3em} and \ \sl L. Vaksman
\thanks{Partially supported by the grant INTAS-94-4720, ISF grant U21200
and grant DKNT-1.4/12}}

\date{\tt Institute for Low Temperarture Physics \& Engineering\\
National Academy of Sciences of Ukraine}

\newtheorem{lemma}{Lemma}
\newtheorem{proposition}[lemma]{Proposition}
\newtheorem{remark}[lemma]{Remark}

\begin{document}

\maketitle
\section{Introduction}

 Consider an irreducible Hermitian symmetric space $X$ of non-compact type.
Let ${\frak g}$ and ${\frak g}_0$ denote the complexifications of the Lie
algebras of the automorphism group of $X$ and the stabilizer of a point $x
\in X$ respectively. Then the center of ${\frak g}_0$ is 1-dimensional
$(Z({\frak g}_0)\;=\:{\Bbb C}\cdot H, \; H \in {\frak g}_0)$, and ${\frak
g}\:=\:{\frak g}_{-1}\oplus{\frak g}_0 \oplus{\frak g}_1$, where ${\frak
g}_{\pm 1}\;=\;\{\xi \in {\frak g}| \:[H,\xi]\,=\,\pm 2 \xi \}$ (see, e.g.,
\cite{He}).

 It was shown by Harish-Chandra that there exists a natural embedding $i:\;X
\hookrightarrow{\frak g}_{-1}$ with $iX$ being a bounded symmetric domain in
${\frak g}_{-1}$ \cite{He}.

 Our purpose is to construct quantum analogues of the (prehomogeneous)
vector space ${\frak g}_{-1}$, the bounded symmetric domain $iX \subset{\frak
g}_{-1}$ and the differential calculus in ${\frak g}_{-1}$.

 Normally we don't dwell on describing the quantum algebras of functions and
quantum exterior algebras in terms of generators and relations, although
that could be done. (The case ${\frak g}\,=\,\frak{sl}_{m+n},\:{\frak
g}_0={\frak s}(\frak{gl}_m \times \frak{gl}_n)$ was partially considered in
\cite{SiV}).

 The simplest homogeneous bounded domain is the unit disc: $U\:=\:\{z
\in{\Bbb C}|\,{\scriptstyle \mid} z {\scriptstyle \mid}\,<\,1\}$. It was
shown in \cite{KL, KRR} that the Poisson brackets \{.,.\} that agree with
the action of the Poisson-Lie group $SU(1,1)$ on $U$ are given by
$$\{z,\overline{z}\}\;=\;i(1-|z|^2)(a+b|z|^2),\quad a,b \in{\Bbb R}.$$

 Our construction (see Section 9) provides a quantization of this bracket
with $b=0$. This "simplest" quantum disc was studied in \cite{NN, VSh}.

 Most of the constructions of this paper originate from the works of V. G.
Drinfeld \cite{Dr1}, S.~Z.~Levendorski\v{i} \& Ya. S. Soibelman \cite{LS}.
Specifically, we follow \cite{Dr1} in replacing the construction of algebras
by forming the dual coalgebras; also our choice of a Poisson cobracket,
together with the associated quantization procedure, is due to \cite{LS}.

 The authors are grateful to V. Akulov, V. Lyubashenko, G. Maltsiniotis,
and D. Shklyarov for the helpful discussion of the results.

\bigskip
\section{Prehomogeneous vector spaces of a commutative
parabolic type}
 Everywhere in the sequel ${\Bbb C}$ will be the ground field.

 Let ${\frak g}$ be a simple complex Lie algebra, ${\frak h}$ its Cartan
subalgebra and $\alpha_i \in {\frak h}^*, \; i=1,\ldots,l$ a simple root
system of ${\frak g}$.

 Choose an element $\alpha_0 \in \{\alpha_i \}_{i=1,\ldots,l}$ and consider
the associated ${\Bbb Z}$-grading
$${\frak g}\;=\;\scriptstyle \bigoplus \limits_{\scriptstyle
j}\textstyle{\frak g}_j,\quad {\frak g}_j\;=\;\{\xi \in{\frak
g}|\:[H_0,\xi]\,=\,2j \xi \},$$
where $H_0 \in {\frak h},\quad \alpha_0(H_0) \,=\,2,\quad
\alpha_i(H_0)\,=\,0$ \ for \ $\alpha_i \,\ne \alpha_0$.

 A subspace ${\frak g}_{-1}$ is called a prehomogeneous vector space of a
commutative parabolic type if the above ${\Bbb Z}$-grading breaks off:
$${\frak g}\;=\;{\frak g}_{-1}\oplus{\frak g}_0 \oplus{\frak
g}_1.\eqno(2.1)$$

 The motives that justify this definition and the list of simple roots
$\alpha \in \{\alpha_i\}_{i=1,\ldots,l}$ with (2.1) being valid are given in
\cite{BR, R}.

 It is worthwhile to note that all the simple roots of series $A_n$ Lie
algebras possess the above property, and for the Lie algebra series
$B_n,C_n,D_n$, together with the exceptional Lie algebras $E_6,E_7$ the set
of such roots is non-void.

 Set ${\frak p}_+\,:=\,{\frak g}_0 \oplus{\frak g}_1,\;{\frak
p}_-\,:=\,{\frak g}_0 \oplus{\frak g}_{-1}$. Our purpose is to construct a
quantum analogue of the graded polynomial algebra ${\Bbb C}[{\frak g}_{-1}]$
on the prehomogeneous vector space ${\frak g}_{-1}$. For this, it would be
useful to have a definition of ${\Bbb C}[{\frak g}_{-1}]$ in terms of the
enveloping algebras $U{\frak g}\supset U{\frak p}_+\supset U{\frak g}_0$
(but not the Lie algebras themselves).

 We start with constructing the coalgebra $V_-$ dual to ${\Bbb C}[{\frak
g}_{-1}]$. Consider the $U{\frak g}$-module $V_-$ determined by its generator
$v \in V_-$ and the relations
$$\xi v_-\,=\,\varepsilon(\xi)v_-,\quad \xi \in U{\frak p}_+\eqno(2.2)$$
where $\varepsilon:U{\frak p}_+\to{\Bbb C}\simeq{\rm End}({\Bbb C})$ is
the trivial representation of $U{\frak p}_+$. Equip $V_-$ with a structure
of a coalgebra \cite{CP} by extending the map $\Delta_-:\,v_-\mapsto
v_-\otimes v_-$ to a morphism of $U{\frak g}$-modules. The existence and
uniqueness of this extension are obvious, and the coassociativity of
$\Delta_-$ follows from
$$(\Delta_-\otimes{\rm id})\Delta_-v_-\,=\,(v_-\otimes v_-)\otimes v_-;\quad
({\rm id}\otimes \Delta_-)\Delta_-v_-\,=\,v_-\otimes(v_-\otimes v_-).$$

 It is easy to verify that $V_-\,=\,\scriptstyle \bigoplus
\limits_{\scriptstyle j}\textstyle(V_-)_j$ with $(V_-)_j\;=\;\{v \in
V_-|\:H_0v\,=\,2jv\}$, and that the dual algebra $\scriptstyle \bigoplus
\limits_{\scriptstyle j}\textstyle((V_-)_j)^*$ to the coalgebra $V_-$ is
canonically isomorphic to ${\Bbb C}[{\frak g}_{-1}]$.

 A replacement of '--' by '+' in the above construction leads to the algebra
of antiholomorphic polynomials on ${\frak g}_{-1}$, which will be denoted by
${\Bbb C}[\overline{\frak g}_{-1}]$. We shall see in the sequel that these
constructions can be transferred to the quantum case where they lead to the
"covariant" algebras ${\Bbb C}[{\frak g}_{-1}]_q,\;{\Bbb C}[\overline{\frak
g}_{-1}]_q$.

\bigskip
\section{Quantum universal enveloping algebras and their
"real forms"}

 It is well known \cite{S} that a simple complex Lie algebra ${\frak g}$
admits a description in terms of generators $\{X_i^\pm,H_i\}_{i=1}^{\;l}$ and
relations
$$[H_i,H_i]=0;\quad[H_i,X_j^{\pm}]=\pm a_{ij}X_j^\pm;\quad
[X_i^+,X_j^-]=\delta_{ij}H_i;\quad{\rm ad}_{X_i^\pm}^{1-a_{ij}}(X_j^\pm)=0.
\eqno(3.1)$$
In the above $i,j \in \{1,\ldots,l \}$, and $(a_{ij})$ is the Cartan matrix
of the simple Lie algebra ${\frak g}$, i. e. $a_{ij}=\alpha_i(H_j)$.

 Let $j_0$ be the number of the simple root $\alpha_0$. The relations (2.2)
can be rewritten in the form
$$\begin{array}{ll}X_j^-v_-=H_jv_-=0,& j=1,2,\ldots,l;\\
                   X_j^+v_-=0, & j \ne j_0.\end{array}$$

 Consider the real Lie subalgebra ${\frak g}(\alpha_0) \subset{\frak g}$
generated by the elements
$$\begin{array}{cccl}X_j^+-X_j^-,&
i(X_j^++X_j^-),& iH_j; & j \ne j_0\\ X_{j_0}^+-X_{j_0}^-, &
                     i(X_{j_0}^++X_{j_0}^-), & iH_{j_0},
\end{array}$$
where $i=\sqrt{-1}$. This subalgebra is interesting because it is the Lie
algebra for the automorphism group of the corresponding bounded symmetric
domain in ${\frak g}_{-1}\subset{\frak g}$. We are seeking for the specific
ways to distinguish $U{\frak g}(\alpha_0)$ inside $U{\frak g}$.

 Recall that $U{\frak g}$ is a Hopf algebra \cite{CP} whose comultiplication
$\Delta$, counit $\varepsilon$ and antipode $S$ are given by
$$\Delta(H_i)=H_i \otimes 1+ 1 \otimes H_i,\quad \Delta(X_i^\pm)=X_i^\pm
\otimes 1+1 \otimes X_i^\pm;$$
$$\varepsilon(H_i)=\varepsilon(X_i^\pm)=0;\qquad S(H_i)=-H_i,\quad
S(X_i^\pm)=-X_i^\pm,$$
$i=1,2,\ldots,l$.

 It is easy to verify that
$$U{\frak g}(\alpha_0)\:=\:\{\xi \in U{\frak g}|\,\xi^*=S(\xi)\},$$
with $*$ being the antilinear involution which depends on $\alpha$ and is
determined by its values on generators $X_j^\pm,\,H_j$ as follows:
\footnote{It is implicit that $(\xi \eta)^*=\eta^*\xi^*,\quad \xi,\eta \in
U{\scriptstyle\frak g}.$}
$$\begin{array}{ccl}H_{j_0}^*=H_{j_0}, & (X_{j_0}^\pm)^*=-X_{j_0}^\mp,\\
                    H_j^*=H_j, & (X_j^\pm)^*=X_j^\mp, & j \ne j_0.
\end{array}\eqno(3.2)$$

 The Hopf algebra $U{\frak g}(\alpha_0)$ doesn't survive under quantization;
in the sequel it will be replaced by the pair $(U{\frak g}, ^*)$. Now let us
consider the quantization of this Hopf $*$-algebra.

 We start with V. G. Drinfeld -- M. Jimbo formulae \cite{Dr1} which
determine a Hopf algebra $U_h{\frak g}$ over ${\Bbb C}[[h]]$ complete in
$h$-adic topology (${\Bbb C}[[h]]$ denotes the ring of formal series). First
of all, choose an invariant scalar product in ${\frak g}$ in such a way that
$d_i=(\alpha_i,\alpha_i)/2>0$. Now $\{X_j^\pm,H_j\}_{j=1,\ldots,l}$ work as
generators of the topological algebra $U_h{\frak g}$, and the list of
relations is as follows:
$$[H_i,H_j]=0,\quad[H_i,X_j^\pm]=\pm a_{ij}X_j^\pm,\quad
[X_i^+,X_j^-]=\delta_{ij}\frac{{\rm sh}(d_jhH_j/2)}{{\rm sh}(d_jh/2)},$$
$$\sum_{k=0}^{1-a_{ij}}(-1)^k\left[\begin{array}{c}1-a_{ij}\\k\end{array}
\right]_h(X_i^\pm)^kX_j^\pm (X_i^\pm)^{(1-a_{ij}-k)}=0.$$

 Here we use the notation
$$\left[\begin{array}{c}n \\ m \end{array}\right]_h \;=\;\prod_{k=1}^n
\frac{{\rm sh}(kh/2)}{{\rm sh}(h/2)}\hbox{\Huge/}\left(\prod_{k=1}^m
\frac{{\rm sh}(kh/2)}{{\rm sh}(h/2)}\cdot \prod_{k=1}^{n-m}
\frac{{\rm sh}(kh/2)}{{\rm sh}(h/2)}\right),$$
$i,j=1,\ldots,l$.

 Comultiplication $\Delta$, counit $\varepsilon$ and antipode $S$ are
determined by their values on the generators:
$$\Delta(H_i)=H_i \otimes 1+1 \otimes H_i;\quad \Delta(X_i^\pm)=X_i^\pm
\otimes e^{hH_id_i/4}+e^{-hH_id_i/4}\otimes X_i^\pm;$$
$$\varepsilon(H_i)=\varepsilon(X_i^\pm)=0;\quad S(H_i)=-H_i;\quad
S(X_i^\pm)=-e^{\pm hd_i/2}\cdot X_i^\pm.$$

 An involution in ${\Bbb C}[[h]]$ is introduced by setting $h^*=h$. We equip
$U_h{\frak g}$ with the structure of $*$-algebra over ${\Bbb C}[[h]]$
defined by (3.2). The pair $(U_h{\frak g},^*)$ will be denoted by $U_h$ for
the sake of brevity.

 A procedure of transition from algebras over ${\Bbb C}[[h]]$ to algebras
over ${\Bbb C}$ is described in \cite{CP}; it allows one to "fix the value
of the formal parameter $h$". Here we only remind the formulae which
describe the "change of variables" corresponding to the generators of the
above algebra:
$$q=e^{-h/2}, \quad K_i^{\pm 1}=e^{\mp hd_iH_i/2},\quad
E_i=X_i^+e^{-hd_iH_i/4}, \quad F_i=e^{hd_iH_i/4}X_i^-.$$

 We fix in the sequel the value of $q \in (0,1)$. The Hopf algebra over
${\Bbb C}$ given by the generators $\{E_i,F_i,K_i^{\pm 1}\}_{i=1}^l$ and the
relations deduced above from the relations in $U_h$, will be denoted by
$U_q{\frak g}$, and the Hopf $*$-algebra $(U_q{\frak g},^*)$ by
$U_q.$\footnote{See the definition of a Hopf $*$-algebra in \cite{CP}.}

 The defining relations for $U_q$ are similar to (3.1), (3.2). We list a
part of them (the quantum analogue of the last among the relations (3.2) can
be found in \cite{L}):
$$K_iK_j=K_jK_i; \quad K_iK_i^{-1}=K_i^{-1}K_i=1; \quad
K_iE_j=q^{d_ia_{ij}}\cdot E_jK_i; \quad K_iF_j=q^{-d_ia_{ij}}F_jK_i;$$
$$E_iF_j-F_jE_i=\delta_{ij}\frac{K_i-K_i^{-1}}{q^{d_i}-q^{-d_i}}$$
$$\Delta(E_i)=E_i \otimes 1+K_i \otimes E_i; \quad \Delta(F_i)=F_i \otimes
K_i^{-1}+1 \otimes F_i; \quad \Delta(K_i)=K_i \otimes K_i;$$
$$\varepsilon(E_i)=\varepsilon(F_i)=\varepsilon(K_i-1)=0;$$
$$S(E_i)=-K_i^{-1}E_i; \quad S(F_i)=-F_iK_i; \quad S(K_i)=K_i^{-1};$$
$$E_j^*=\left\{\begin{array}{r|c}K_jF_j & j \ne j_0 \\ -K_jF_j & j=j_0
\end{array}\right., \qquad F_j^*=\left\{\begin{array}{r|c}E_jK_j^{-1} & j
\ne j_0 \\ -E_jK_j^{-1} & j=j_0 \end{array}\right.,$$
$K_j^*=K_j,\quad i,j \in \{1,\ldots,l\}$.

 We equip the Hopf algebra $U_q{\frak g}$ with a grading as follows:
$${\rm deg}\,K_j\:=\:{\rm deg}\,E_j\:=\:{\rm deg}\,F_j\:=\:0,\qquad j \ne
j_0$$
$${\rm deg}\,K_{j_0}\,=\,0,\quad {\rm deg}\,E_{j_0}\,=\,1,\quad {\rm
deg}\,F_{j_0}\,=\,-1.$$

\bigskip
\section{Covariant algebras and involutions}

 Remind that ${\Bbb C}$ is endowed with a structure of a $U_q{\frak
g}$-module by means of a counit $\varepsilon:\,U_q{\frak g}\,\to \,{\Bbb
C}\simeq {\rm End}({\Bbb C})$.

 Let ${\cal F}$ be a unital algebra over ${\Bbb C}$, which is also a
$U_q{\frak g}$-module. We call ${\cal F}$ a {\sl $U_q{\frak g}$-module
algebra} if the multiplication
$$m:\,{\cal F} \otimes {\cal F} \to {\cal F};\quad m:\,f_1 \otimes f_2
\mapsto f_1f_2,\quad f_1,f_2 \in {\cal F}$$
and the unit
$$1:\,{\Bbb C}\to {\cal F};\quad z \mapsto z \cdot 1, \quad z \in {\Bbb C}$$
are morphisms of $U_q{\frak g}$-modules.\footnote{${\cal F} \otimes {\cal
F}$ becomes a $U_q{\scriptstyle \frak g}$-module by setting $\xi(f_1 \otimes
f_2)\,=\,\sum \limits_j \xi_j^\prime f_1 \otimes \xi_j^{\prime \prime}f_2$
for $\xi \in U_q{\scriptstyle \frak g}$ with $\Delta(\xi)\,=\,\sum \limits_j
\xi_j^\prime \otimes \xi_j^{\prime \prime},\;f_1,f_2 \in {\cal F}$.}

 Together with the term "$U_q{\frak g}$-module algebra" we shall elaborate
the substitute term {\sl "covariant algebra"} for the sake of brevity in the
cases when no confusion can occur.

 Covariant modules and covariant bimodules over covariant algebras are
defined in a similar way (see \cite{A, W}).

 An involutive $({\cal F},*)$ algebra is said to be covariant \cite{SoV} if
it is a $U_q{\frak g}$-module algebra and for all $\xi \in U_q{\frak g}, \,
f \in {\cal F}$ one has
$$(\xi f)^*\,=\,(S(\xi))^*f^*.\eqno(4.1)$$

 A linear functional $\nu:\,{\cal F} \to {\Bbb C}$ is called an invariant
integral if
$$\nu(\xi f)\,=\,\varepsilon(\xi)\nu(f), \quad \xi \in U_q{\frak g},\;f \in
{\cal F}.$$

 The "compatibility condition" for involutions (4.1) is extremely important
since it allows one to use the "positive" invariant integrals for producing
$*$-representations of $U_q{\frak g}$ in the "Hilbert function spaces":
$$(f_1,f_2)\,=\,\nu(f_2^*f_1), \quad f_1,f_2 \in {\cal F}.$$

 The problem of decomposing such $*$-representations is a typical one in
harmonic analysis. On this way, for instance, the Plancherel measure for
quantum $SU(1,1)$ was found (see \cite{SoV}).

\bigskip
\section{Generalized Verma modules}

 Choose a linear functional $\lambda \in{\frak h}^*$ so that
$m_j=\lambda(H_j)$ are non-positive integers for $j \ne j_0$.

 Consider the graded $U_q{\frak g}$-module determined by the single
generator $v_+(\lambda) \in V_+(\lambda)$ and the relations
$$F_iv_+(\lambda)=0,\; K_i^{\pm 1}v_+(\lambda)=e^{\mp
d_im_ih/2}v_+(\lambda),\quad i=1,\ldots,l;$$
$$E_j^{-m_j+1}v_+(\lambda)=0,\quad j \ne j_0;$$
$${\rm deg}(v_+(\lambda))={1 \over 2}\lambda(H_0).$$

 Note that $V_+(\lambda)\,=\,\scriptstyle \bigoplus \limits_{\scriptstyle j}
\textstyle V_+(\lambda)_j$, with $V_+(\lambda)_j\:=\:\{v \in
V_+(\lambda)|\,{\rm deg}(v)=j \}$, and ${\rm dim}V_+(\lambda)_j<\infty$.

 The finite dimensionality of the homogeneous component $V_+(\lambda)_j$
follows from the decomposition
$$V_+(\lambda)_j \:=\:\bigoplus_{\{\mu \in{\frak
h}^*|\,\mu(H_0)=2j\}}V_+(\lambda)_\mu$$
into a finite sum of the finite dimensional weight subspaces
$$V_+(\lambda)_{\mu}\:=\:\{v \in V_+(\lambda)|\,K_jv=e^{-d_j \mu(H_j)h/2}v,\,
j=1,\ldots,l\}.$$

 The graded modules $V_-(\lambda)$ are defined in a similar way:
$$E_iv_-(\lambda)=0,\quad K_i^{\pm
1}v_-(\lambda)=e^{\mp d_im_ih/2}v_-(\lambda),\quad i=1,\ldots,l$$
$$F_j^{m_j+1}v_-(\lambda)=0,\quad j \ne j_0;\qquad {\rm deg}(v_-(\lambda))={1
\over 2}\lambda(H_0).$$

 Now suppose $m_{j_0}=\lambda(H_{j_0}) \in {\Bbb Z}$.

 Consider the longest element $w_0$ of the Weyl group $W$ for a Lie
algebra ${\frak g}$. It is very well known from \cite{CK, CP} that one can
associate to each reduced decomposition of $w_0$ a Poincar\`e-Birkhoff-Witt
basis in $U_q{\frak g}$. We demonstrate the reduced decompositions for which
this basis "generates" the bases of weight vectors in generalized Verma
modules.

 Let ${\frak g}^\prime \in {\frak g}$ be a Lie subalgebra generated by
$\{X_j^\pm,H_j\}_{j \ne j_0}$, and let also $W^\prime \in W$ be a subgroup
generated by simple reflections $s(\alpha_j),\: j \ne j_0$. Obviously,
$W^\prime$ is a Weyl group of the Lie algebra ${\frak g}^\prime$.

 Denote by $U \in W$ the subset of such elements $u \in W$ that
$$l(s(\alpha_j)u) > l(u) \quad \hbox{for all} \quad j \ne j_0.$$

 It is known from \cite[p. 19]{Hu}, that, firstly, each element $w
\in W$ admits the unique decomposition $w=w^\prime \cdot u$ with $w^\prime
\in W^\prime,\;u \in U$. Secondly, if $w^\prime \in W^\prime,\; u \in U$,
then one has
$$l(w^\prime \cdot u)\,=\,l^\prime(w^\prime)+l(u),$$
with $l^\prime(w^\prime)$ being the length of the element $w^\prime$ in
$W^\prime$, and $l(u),\; l(w^\prime u)$ the lengths of $u,\,w^\prime u$ in
$W$.

 That is, one can find in $U$ the unique element $u_0$ of maximum length
such that $w_0=w_0^\prime \cdot u_0$. ($w_0^\prime$ here is the longest
element of $W^\prime$.) Now one can derive the desired reduced
decompositions of $w_0$ by multiplication from the reduced decompositions
of $w_0^\prime$ and $u$.

\bigskip
\section{From coalgebras to algebras}

 Let $U_q{\frak g}^{op}$ stand for the Hopf algebra derived from $U_q{\frak
g}$ by replacing its comultiplication by the opposite one.

 We intend to use the generalized Verma modules for producing coalgebras
dual to covariant algebras. To provide a precise correspondence between
these two notions, we are going to replace $U_q{\frak g}$ by $U_q{\frak
g}^{op}$ in tensor products of generalized Verma modules.

Consider the $U_q{\frak g}$-modules $V_ \pm(0)$. Evidently the maps
$$\Delta_ \pm:\:v_ \pm(0)\mapsto v_ \pm(0)\otimes v_ \pm(0);\quad
\varepsilon_ \pm:\: v_ \pm(0)\mapsto 1$$
admit the unique extensions to morphisms of $U_q{\frak g}$-modules:
$$\Delta_ \pm:\:V_ \pm(0)\,\to \,V_ \pm(0)\otimes V_ \pm(0);\quad
\varepsilon_ \pm:\:V_ \pm(0)\to{\Bbb C}.$$

 Just as in the case $q=1$, one can verify that the operations $\Delta_ \pm$
are coassociative, and that $\varepsilon_ \pm$ are the counits for
coalgebras respectively with $\Delta_ \pm$.

 It follows that the vector spaces $(V_
\pm(0))^*\:\stackrel{def}{=}\:\scriptstyle \bigoplus \limits_{\scriptstyle
j}\textstyle (V_ \pm(0)_j)^*$ are covariant algebras.\footnote{Remind
\cite{CP} that the dual $U_q{\scriptstyle \frak g}$-module structure is
given by $\xi f(v)\:\stackrel{def}{=}\:f(S(\xi)v)$, with $\xi \in
U_q{\scriptstyle \frak g},\;v \in V_ \pm(0),\; f \in V_ \pm(0)^*$.}
Introduce the notation
$${\Bbb C}[{\frak g}_{-1}]_q \:=\:V_-(0)^*,\quad {\Bbb C}[\overline{\frak
g}_{-1}]_q \:=\:V_+(0)^*.$$
These covariant algebras may be treated as q-analogues of polynomial
algebras (holomorphic or antiholomorphic identified by the sign) on the
quantum prehomogeneous space ${\frak g}_{-1}$.

\bigskip
\section{Polynomial algebra}

 Consider the algebra ${\rm Pol}({\frak g}_{-1})={\Bbb C}[{\frak
g}_{-1}]\otimes{\Bbb C}[\overline{\frak g}_{-1}]$ of all polynomials on
${\frak g}_{-1}$. Holomorphic and antiholomorphic polynomials admit the
embeddings into this algebra as follows:

$${\Bbb C}[{\frak g}_{-1}]\hookrightarrow {\Bbb C}[{\frak g}_{-1}]\otimes
{\Bbb C}[\overline{\frak g}_{-1}],\quad f \mapsto f \otimes 1,$$
$${\Bbb C}[\overline{\frak g}_{-1}]\hookrightarrow {\Bbb C}[{\frak
g}_{-1}]\otimes {\Bbb C}[\overline{\frak g}_{-1}],\quad f \mapsto 1 \otimes
f.$$

 Our desire is to obtain that sort of algebra and similar embeddings in the
quantum case $(q \ne 1)$. For that, we intend to equip the $U_q{\frak
g}$-module ${\rm Pol}({\frak g}_{-1})_q \:\stackrel{def}{=}\:{\Bbb C}[{\frak
g}_{-1}]_q \otimes{\Bbb C}[\overline{\frak g}_{-1}]_q$ with a structure of
covariant algebra in such a way that the maps $f \mapsto f \otimes 1,\; f
\mapsto 1 \otimes f$ turn out to be algebra homomorphisms.

 Our approach is completely standard \cite{JS}. Define the product of
$\varphi_+ \otimes \varphi_-,\,\psi_+ \otimes
\psi_-\:\in\:{\rm Pol}({\frak g}_{-1})_q$ as follows:
$$(\varphi_+ \otimes \varphi_-)(\psi_+ \otimes \psi_-)\,=\,m_+ \otimes
m_-(\varphi_+ \otimes \check{R}(\varphi_- \otimes \psi_+)\otimes \psi_-).$$
Here $m_+,\,m_-$ are the multiplications in ${\Bbb C}[{\frak
g}_{-1}]_q,\,{\Bbb C}[\overline{\frak g}_{-1}]_q$ respectively, and
$\check{R}:\;{\Bbb C}[\overline{\frak g}_{-1}]_q \otimes {\Bbb C}[{\frak
g}_{-1}]_q \to {\Bbb C}[{\frak g}_{-1}]_q \otimes {\Bbb C}[\overline{\frak
g}_{-1}]_q$ is the morphism of $U_q{\frak g}$-modules defined below by V. G.
Drinfeld's universal R-matrix \cite{Dr1}.

  One can find in \cite{Dr1, Dr2} the description of properties of the
universal R-matrix which unambiguously determine it as an element of an
appropriate completion of $U_h{\frak g}\otimes U_h{\frak g}$. In particular,
$$S \otimes S(R)=R, \quad R^{*\otimes *}=R^{21}.\eqno(7.1)$$
The latter relation involves the element $R^{21}$ which is derived from $R$
by permutation of tensor multiples. The proof of this relation is completely
similar to that of Proposition 4.2 in \cite{Dr2}.

 In \cite{CP} there is an explicit "multiplicative" formula for the
universal R-matrix. More precisely, any reduced decomposition of the maximum
length element $w_0$ possesses its own multiplicative formula. In the sequel
we intend to restrict ourselves to those reduced decompositions which come
from Section 5. (Note that the "multiplicative" formula was discovered in
the papers of S. Z. Levendorski\v{i}, Ya. S. Soibelman and also by A. N.
Kirillov, N. Yu. Reshetikhin, see \cite{CP}. Its application should take
into account the inessential differences in the choice of generators and
deformation parameters in this work as compared with \cite{CP}.
Specifically, one has to substitute $X_i^+,\,X_i^-,\,H_i,\,h,\,K_i,\,q$ \ by
\ $-S(E_i),\,-S(F_i),\,-S(H_i),\,h/2,\,K_i^{-1},\,q^{-1}$.)

 It is easy to show that the universal R-matrix determines a linear operator
in ${\Bbb C}[\overline{\frak g}_{-1}]_q \otimes {\Bbb C}[{\frak g}_{-1}]_q$.

 Now we are in a position to define the operator $\check{R}$ in a standard
way: $\check{R}=\sigma \cdot R$ with $\sigma:\:a \otimes b \mapsto b \otimes
a$ being a permutation of tensor multiples. Thus $\check{R}$ becomes a
morphism of $U_q{\frak g}$-modules since \cite{Dr1, Dr2, CP}
$$\Delta^{op}(\xi)\,=\,R \Delta(\xi)R^{-1}, \quad \xi \in U_h{\frak g}.$$

 The associativity of the multiplication in ${\rm Pol}({\frak g}_{-1})_q$
can be easily derived by the standard argument \cite{JS, K} from the
relations
$$(\Delta \otimes {\rm id})(R)\,=\,R^{13}R^{23}, \quad ({\rm id}
\otimes\Delta)(R)\,=\,R^{13}R^{12}.$$
(Here $R^{12}\,=\,\sum \limits_i a_i \otimes b_i \otimes 1,\;
R^{23}\,=\,\sum \limits_i 1 \otimes a_i \otimes b_i,\; R^{13}\,=\,\sum
\limits_i a_i \otimes 1 \otimes b_i$ whenever $R\,=\,\sum \limits_i a_i
\otimes b_i$, see \cite{Dr1, Dr2, CP}).

 The existence of a unit and covariance of ${\rm Pol}({\frak g}_{-1})_q$ are
evident.

\bigskip
\section{Involution}

 Consider the antilinear operators $^*:V_+(0) \to V_-(0);\quad ^*:V_-(0) \to
V_+(0)$, which are determined by their properties as follows. Firstly, $v_
\pm(0)^*\,=\,v_ \mp(0)$ and, secondly,
$$(\xi v)^*\,=\,(S^{-1}(\xi))^*v^* \eqno(8.1)$$
for all $v \in V_ \pm(0),\; \xi \in U_q{\frak g}$.

 To rephrase the above, we set up
$$(\xi v_ \pm(0))^*\,=\,(S^{-1}(\xi))^*v_ \pm(0)^*.$$
It follows from the definition of $V_ \pm(0)$ that the involution as above
is well defined. In particular, (8.1) can be easily deduced; it also follows
from the relation $(S^{-1}((S^{-1}(\xi))^*))^*\,=\,\xi$ that the operators
constructed above are mutually converse.

 The duality argument allows one to form the mutually converse
antihomomorphisms $^*:{\Bbb C}[{\frak g}_{-1}]_q \to {\Bbb C}[\overline{\frak
g}_{-1}]_q;\quad ^*:{\Bbb C}[\overline{\frak g}_{-1}]_q \to {\Bbb C}[{\frak
g}_{-1}]_q:$
$$f^*(v) \stackrel{def}{=}\overline{f(v^*)},\quad v \in V_ \pm(0),\:f \in
V_ \pm(0)^*.\eqno(8.2)$$

 Now we are in a position to define the antilinear operator $*$ in ${\rm
Pol}({\frak g}_{-1})_q$ by
$$(f_+ \otimes f_-)^* \stackrel{def}{=} f_-^* \otimes f_+^*,$$
for $f_+ \in {\Bbb C}[{\frak g}_{-1}]_q,\: f_- \in {\Bbb C}[\overline{\frak g
}_{-1}]_q$, and also to show that it equips ${\rm Pol}({\frak g}_{-1})_q$
with a structure of covariant involutive algebra.

 It remains to verify that $*$ is an antihomomorphism of ${\rm Pol}({\frak
g}_{-1})_q$. The best way to prove the relation
$$(f_1 f_2)^*\,=\,f_2^*f_1^*; \quad f_1,f_2 \in {\rm Pol}({\frak g}_{-1})_q$$
is in applying (7.1) and the duality argument described in details in the
concluding section of the present paper. (Note that it suffices to prove the
relation $(f_1f_2)^*(v)\,=\,f_2^*f_1^*(v)$ for the generator $v=v_-(0)
\otimes v_+(0)$ of the $U_q{\frak g}$-module $V_-(0) \otimes V_+(0)$ since
the map $\xi \mapsto (S^{-1}(\xi))^*$ is an antiautomorphism of the
coalgebra $U_q{\frak g}$).

 Verify (4.1); it suffices to consider the case $f \in V_ \pm(0)^*$. An
application of (8.2) and the relation $S((S(\xi))^*)\,=\,\xi^*,\;\xi \in
U_q{\frak g}$, yields for all $v \in V_ \pm(0),\; f \in V_ \pm(0)^*$
$$f((\xi v)^*)\,=\,f((S^{-1}(\xi))^*v^*)$$
$$\overline{f(S(\xi)v^*)}\,=\,\overline{f((\xi^*v)^*)}$$
$$\overline{(\xi f)(v^*)}\,=\,f^*(\xi^*v)$$
$$\overline{(\xi f)(v^*)}\,=\,f^*(S((S(\xi))^*)v)$$
$$(\xi f)^*(v)\,=\,((S(\xi))^*f^*)(v).$$
Thus, in the special case $f \in V_ \pm(0)^*$ (4.1) is proved. Hence it is
also valid for all $f \in {\rm Pol}({\frak g}_{-1})_q$ since the antipode is
an antiautomorphism of the coalgebra $U_q{\frak g}$ and the involution $*$
is its automorphism. In fact, if $f=f_+f_-,\; f_ \pm \in(V_ \pm(0))^*$ and
$\Delta(\xi)\,=\,\sum \limits_j \xi_j^\prime \otimes \xi_j^{\prime
\prime};\; \xi_j^\prime,\xi_j^{\prime \prime} \in U_q{\frak g}$ then one has
$$(\xi(f_+f_-))^*\,=\,\sum_j(\xi_j^{\prime \prime}f_-)^*(\xi_j^\prime
f_+)^*,$$
$$(S(\xi))^*(f_+f_-)^*\,=\,(S(\xi))^*(f_-^*f_+^*)\,=\,
\sum_j((S(\xi_j^{\prime \prime}))^*f_-^*)((S(\xi_j^\prime))^*f_+^*).$$

\bigskip
\section{The simplest example}

 Let ${\frak g}=\frak{sl}_2$, then one has ${\frak g}={\frak
g}_{-1}\oplus{\frak g}_0 \oplus{\frak g}_1$, with ${\frak g}_0$ and ${\frak
g}_{\pm 1}$ being Cartan and Borel subalgebras of $\frak{sl}_2$
respectively. In particular, ${\rm deg}({\frak g}_{-1})=1$.

 The algebra $U_q \frak{sl}_2$ is given by its generators $K^{\pm
1},\,E,\,F$ and the relations
$$KK^{-1}=K^{-1}K=1,\quad K^{\pm 1}E=q^{\pm 2}EK^{\pm 1},\quad K^{\pm
1}F=q^{\mp 2}FK^{\pm 1},$$
$$EF-FE=(K-K^{-1})/(q-q^{-1}).$$

 Remind that comultiplication $\Delta$, counit $\varepsilon$ and antipode $S$
are defined on the above generators as follows:
$$\Delta(E)=E \otimes 1\,+\,K \otimes E,\quad \Delta(F)=F \otimes
K^{-1}\,+\,1 \otimes F,\quad \Delta(K^{\pm 1})=K^{\pm 1} \otimes K^{\pm 1};$$
$$\varepsilon(E)=\varepsilon(F)=0,\quad \varepsilon(K^{\pm 1})=1;$$
$$S(E)=-K^{-1}E,\quad S(F)=-FK,\quad S(K^{\pm})=K^{\mp}.$$

 In the notation
$$q=e^{-h/2},\quad K^{\pm 1}=e^{\mp hH/2},\quad E=X^+e^{-hH/4},\quad
F=e^{hH/4}X^-$$
the V. G. Drinfeld's formula for the universal R-matrix \cite{Dr1} acquires
the form
$$R={\rm exp}_{q^2}((q^{-1}-q)E \otimes F)\cdot {\rm exp}(H \otimes H \cdot
h/4)$$
with $\exp_t(x)\:=\:\sum \limits_{n=0}^\infty x^n(\prod \limits_{j=1}^n
\frac{\textstyle 1-t^j}{\textstyle 1-t})^{-1}$.

 The involution $*$ in $U_q \frak{su}(1,1)=(U_q \frak{sl}_2,^*)$ is defined
on the generators $E,\,F,\,K^{\pm 1}$ by
$$E^*=-KF,\quad F^*=-EK^{-1},\quad (K^{\pm 1})^*=K^{\pm 1}$$
(equivalently, on the generators $X^\pm,\,H$ of $U_h \frak{sl}_2$ it is
defined by $(X^\pm)^*=-X^\mp,\; H^*=H$).

 Consider the $U_q \frak{sl}_2$-module $V_+(0)$ determined by its single
generator $v_+(0) \in V_+(0)$ and the relations $Fv_+(0)=0,\;K^{\pm
1}v_+(0)=v_+(0)$. This module admits the decomposition
$$V_+(0)\,=\,\scriptstyle \bigoplus \limits_{\scriptstyle j \in{\Bbb Z}_+}
\textstyle V_+(0)_j,\quad V_+(0)_j={\Bbb C}\cdot E^j \cdot v_+(0).$$
Hence $\{E^jv_+(0)\}_{j \in {\Bbb Z}_+}$ is a basis in $V_+(0)$.

 Define a linear functional $a_- \in V_+(0)^*={\Bbb C}[\overline{\frak
g}_{-1}]_q$ by \\
$a_-(S(E^j)v_+(0))\,=\,\left\{\begin{array}{c|c}1 & j=1 \\ 0 & j \ne 1
\end{array} \right.$.

 Prove that $a_-$ is a generator of ${\Bbb C}[\overline{\frak g}_{-1}]_q$,
and that for any polynomial $P \in{\Bbb C}[t]$
$$K^{\pm 1}:\:P(a_-)\mapsto P(q^{\mp 2}a_-),\eqno(9.1)$$
$$E:\:P(a_-)\mapsto(D_-P)(a_-)\eqno(9.2)$$ $$F:\:P(a_-)\mapsto -q \cdot
a_-^2 \cdot (D_+P)(a_-)\eqno(9.3)$$
where $(D_ \pm P)(t)=(P(q^{\pm 2}t)-P(t))/(q^{\pm 2}t-t)$.

 Note first that the relations
$$Ka_-(S(E^j)v_+(0))=\left\{\begin{array}{c|c}q^{-2} & j=1 \\ 0 & j \ne 1
\end{array}\right.,$$
$$Ea_-(S(E^j)v_+(0))=\left\{\begin{array}{c|c}1 & j=0 \\ 0 & j \ne 0
\end{array}\right.$$
imply that
$$Ka_-=q^{-2}a_-,\quad Ea_-=1.\eqno(9.4)$$

 Now apply the covariance of ${\Bbb C}[\overline{\frak g}_{-1}]_q$ to obtain
$$K^{\pm 1}(P_1(a_-)P_2(a_-))=K^{\pm 1}(P_1(a_-))\cdot K^{\pm 1}(P_2(a_-)),$$
$$E(P_1(a_-)P_2(a_-))=E(P_1(a_-))\cdot P_2(a_-)+K(P_1(a_-))\cdot
E(P_2(a_-))$$
for any polynomials $P_1,\,P_2$. This already allows one to deduce (9.1)
(9.2) from (9.4).

 It is worthwhile to note that $a_-^j \ne 0$ for all $j \in {\Bbb Z}_+$
since
$$E^ja_-^j=\prod_{k=1}^j((q^{-2k}-1)/(q^{-2}-1)) \ne 0.$$
This implies that $(V_+(0)_j)^*={\Bbb C}\cdot a_-^j$. Hence $\{a_-^j\}_{j
\in {\Bbb Z}_+}$ is a basis of the vector space ${\Bbb C}[\overline{\frak
g}_{-1}]_q$. That is, $a_-$ is a generator of the algebra ${\Bbb
C}[\overline{\frak g}_{-1}]_q$.

 Now prove (9.3) in the special case $P(a_-)=a_-$. Specifically, we are
going to demonstrate
$$Fa_-=-q \cdot a_-^2.\eqno(9.5)$$
Since $Fa_- \in(V_+(0)_2)^*={\Bbb C}a_-^2$ we have $Fa_-={\rm const}\cdot
a_-^2$. The fact that the constant in the latter relation is $-q$ follows
easily from
$$E(Fa_-)\,=\,\frac{K-K^{-1}}{q-q^{-1}}a_-\,=\,-(q+q^{-1})a_-,$$
$$E(a_-^2)\,=\,(q^{-2}+1)a_-\,=\,q^{-1}(q^	{-1}+q)a_-$$
together with $a_- \ne 0$.

 The passage from the special case $P(a_-)=a_-$ to the general case can be
performed (just as above) by a virtue of covariance. Specifically,
$$F(P_1P_2)\:=\:F(P_1)\cdot K^{-1}(P_2)\,+\,P_1F(P_2)$$
for any "polynomials" $P_1(a_-),\,P_2(a_-)$.

 Now turn to the description of the covariant algebra ${\Bbb C}[{\frak
g}_{-1}]_q$ for the same case ${\frak g}=\frak{sl}_2$. One has
$V_-(0)\:=\:\scriptstyle \bigoplus \limits_{\scriptstyle -j \in{\Bbb
Z}_+}\textstyle V_-(0)_j,\; V_-(0)_{-j}\,=\,{\Bbb C}\cdot F^jv_-(0)$.
Define also the "coordinate function" $a_+$ by
$$a_+(S(F^j)v_-(0))\,=\,\left\{\begin{array}{c|c}1 & j=1 \\ 0 & j \ne 1
\end{array}\right..$$

 Now one can prove in the same way as above that $a_+$ is the generator of
${\Bbb C}[{\frak g}_{-1}]_q$ and
$$K^{\pm 1}:\:P(a_+)\mapsto P(q^{\pm 2}a_+),$$
$$F:\:P(a_+)\mapsto(D_-P)(a_+),$$
$$E:\:P(a_+)\mapsto -qa_+^2\cdot(D_+P)(a_+)$$
for any polynomial P of a single indeterminate.

 In particular, one has
$$K^{\pm 1}a_+=q^{\pm 2}a_-,\quad Fa_+=1,\quad Ea_+=-qa_+^2.\eqno(9.6)$$

 Note that if $\{f_i\}$ are the generators of a covariant algebra ${\cal F}$
and $\{a_j\}$ the generators of a Hopf algebra $A$, then the action of $A$
on ${\cal F}$ can be unambiguously retrieved from the action of $\{a_j\}$ on
$\{f_i\}$.

 Turn to the description of the covariant algebra ${\rm Pol}({\frak
g}_{-1})_q$ in terms of generators and relations.

 By our construction, the covariant algebras ${\Bbb C}[{\frak g}_{-1}]_q$ and
${\Bbb C}[\overline{\frak g}_{-1}]_q$ are embedded into ${\rm Pol}({\frak
g}_{-1})_q$.

 It follows from the explicit formula for the universal R-matrix and the
definition of the action of $\exp(H \otimes Hh/4)$ on the weight vectors
that
$$e^{H \otimes Hh/4}a_- \otimes a_+\:=\:q^{-{1 \over 2}\cdot 2(-2)}\cdot
a_- \otimes a_+,$$
$$a_-a_+\:=\:q^2(a_+a_-\,+\,q^{-1}(1-q^2)Fa_+\cdot Ea_-).$$
Finally we have:
$$a_-a_+\:=\:q^2a_+a_-\,+\,q(1-q^2).\eqno(9.7)$$

 Since ${\rm Pol}({\frak g}_{-1})_q \:=\:{\Bbb C}[{\frak g}_{-1}]_q \otimes
{\Bbb C}[\overline{\frak g}_{-1}]_q$, we deduce that (9.7) gives a complete
list of relations between the generators $\{a_+,a_-\}$ of ${\rm Pol}({\frak
g}_{-1})_q$, that is the natural map ${\Bbb C}\langle
a_+,a_-\rangle/(a_-a_+-(q^2a_+a_-+q(1-q^2)))\to{\rm Pol}({\frak g}_{-1})_q$
is injective. The action of the generators $\{K^{\pm 1},E,F\}$ of $U_q{\frak
g}$ on the generators $\{a_+,a_-\}$ of ${\rm Pol}({\frak g}_{-1})_q$ is given
by (9.4) -- (9.6).

 It remains to describe the involution $*$.

 We start with proving that
$$a_+^*={\rm const}\cdot a_-,\quad a_-^*={\rm const}\cdot a_+,\eqno(9.8)$$
and then find the constants by comparing the explicit expressions for
$Ea_+^*$ and $Ea_-$. The relations (9.8) follow from the decompositions
${\Bbb C}[{\frak g}_{-1}]_q\:=\:\scriptstyle \bigoplus \limits_{\scriptstyle
i}\textstyle(V_-(0)_i)^*, \quad {\Bbb C}[\overline{\frak
g}_{-1}]_q\:=\:\scriptstyle \bigoplus \limits_{\scriptstyle
i}\textstyle(V_+(0)_i)^*$ and $(V_ \mp(0)_{\pm 1})^*\,=\,{\Bbb C}\cdot a_
\pm,\quad ^*:\,(V_ \pm(0)_i)^*\to(V_\mp(0)_{-i})^*$.

 It was pointed out before that $Ea_-=1$. Let's compute $Ea_+^*$. First use
the relation
$$(S(F))^*\,=\,(-FK)^*\,=\,-K^*F^*\,=\,-K \cdot(-EK^{-1})\,=\,q^2E$$
and the compatibility condition (4.1) for involutions to obtain
$$q^2E  a_+^*\,=\,(S(F))^*a_+^*\,=\,(Fa_+)^*\,=\,1^*\,=\,1.$$
Thus we have $q^2 \cdot Ea_+^*\,=\,Ea_-$. Now (9.8) implies
$$a_+^*=q^{-2}a_-;\quad a_-^*=q^2a_+.\eqno(9.9)$$

 The only shortcoming of the definition of the covariant $*$-algebra ${\rm
Pol}({\frak g}_{-1})_q$ is that it is excessively abstract. In the example
for $U_q \frak{su}(1,1)$ we got another description of that covariant
$*$-algebra.  Specifically, its generators are $\{a_+,a_-\}$, its complete
list of relations reduces to (9.7), the action of $U_q \frak{su}(1,1)$ is
given by (9.4) -- (9.7), and the involution is determined by (9.9).

 Note that in the work \cite{VSh} on the function theory in the unit disc the
generator $z=q^{1/2}\cdot a_+$ was implemented instead of $a_+$. In this
setting, (9.9) implies the relation $z^*=q^{-3/2}a_-$, and (9.7) can be
rewritten as
$$z^*z-q^2zz^*\,=\,1-q^2.\eqno(9.10)$$
(The substitution $q=e^{-h/2}$ and the \underline{formal} passage to a limit
as $h \to 0$ yield (cf. (1.1)): $\lim \limits_{h \to
0}\frac{\textstyle[z,z^*]}{\textstyle ih}\,=\,i(1-zz^*)$.)

\bigskip
\section{Quantum disc and other bounded symmetric
domains}

 Proceed with studying the $*$-algebra ${\rm Pol}({\frak g}_{-1})_q$ which
was under investigation in the previous section. Evidently, the formulae
$$T_ \varphi(z)=e^{i \varphi},\quad T_ \varphi(z^*)=e^{-i \varphi},\quad
\varphi \in {\Bbb R}/2 \pi{\Bbb Z}$$
determine the one-dimensional representation of ${\rm Pol}({\frak
g}_{-1})_q$.  We shall also need a faithful infinitely dimensional
$*$-representation $T$ in the Hilbert space $l^2({\Bbb Z}_+)$ given by
$$T(z)e_m \,=\,(1-q^{2(m+1)})^{1/2}e_{m+1},$$
$$T(z^*)e_{m+1}\,=\,(1-q^{2(m+1)})^{1/2}e_m,$$
$$T(z^*)e_0\,=\,0,$$
with $\{e_m\}_{m \in {\Bbb Z}_+}$ being the standard basis in $l^2({\Bbb
Z}_+)$. An application of the standard techniques of operator theory in
Hilbert spaces \cite{VSo} allows one to prove that any irreducible
$*$-representation of the above algebra is unitarily equivalent to one of
the representations $\{T_ \varphi \}_{\varphi \in{\Bbb R}/2 \pi{\Bbb
Z}},\,T$.

 Note that the spectrum of $T(z)$ is the closure $\overline{U}$ of the unit
disc $U$ in ${\Bbb C}$. Just as in \cite{VSo}, we use the notion "algebra of
continuous functions in the quantum disc" for a completion of ${\rm
Pol}({\frak g})_q$ with respect to the norm $\parallel f \parallel\,=\,{\rm
sup}\parallel \rho(f)\parallel$. Here $\rho$ varies inside the class of all
irreducible $*$-representations up to unitary equivalence. One can easily
deduce from the above that $\parallel f \parallel \,=\, \parallel Tf
\parallel$.

 The enveloping von Neumann algebra \cite{Di} of the above $C^*$-algebra
will be denoted by $L^\infty(U)_q$ and called the algebra of continuous
functions in the quantum disc. Certainly, {\sl $L^\infty$ is worthwhile not
alone, but only together with a distinguished dense covariant subalgebra
${\rm Pol}({\frak g}_{-1})_q$} (cf. \cite{W}).

 Note that our quantum disc is only one among those described in \cite{KL}.
Others can be derived from this one by a standard argument normally referred
to as quantization by Berezin \cite{VSh}.

 Remark also that the definition of $L^\infty(U)_q$ which implements a
completion procedure and passage to an enveloping von Neumann algebra
doesn't use the specific features of the special case ${\frak
g}=\frak{sl}_2$. That is, to any irreducible prehomogeneous vector space of
commutative parabolic type we associate a pair constituted by a von Neumann
algebra $L^\infty(U)_q$ and its dense covariant subalgebra ${\rm Pol}({\frak
g}_{-1})_q$.

\bigskip
\section{Differential calculi: the outline}

 We follow G. Maltsiniotis \cite{M} in choosing the basic idea of producing
the differential calculi. Specifically, we first construct differential
calculi of order one, and then embed them into complete differential calculi
by a simple argument described in \cite{M}, the proof of theorem (1.2.3).

 To outline the construction of order one differential calculi, we restrict
ourselves to the simplest example of a quantum prehomogeneous vector space.

 At the first step we consider the type (1,0) forms with holomorphic
coefficients $f \cdot dz,\,f \in {\Bbb C}[{\frak g}_{-1}]_q$, and type (0,1)
forms with antiholomorphic coefficients $f \cdot dz^*,\,f \in {\Bbb
C}[\overline{\frak g}_{-1}]_q$. We prove that
$$dz \cdot z\,=\,q^2z \cdot dz;\quad dz^* \cdot z^*\,=\,q^{-2}z^* \cdot
dz^*.\eqno(11.1)$$

 At the second step we assume the consideration of all the forms of types
(1,0) and (0,1): $fdz,\,fdz^*,\:f \in{\rm Pol}({\frak g}_{-1})_q$. We prove
that $$dz \cdot z^*\,=\,q^{-2}z^*\cdot dz;\quad dz^*\cdot z \,=\,q^2z \cdot
dz^*.\eqno(11.2)$$

 At the third step we turn to higher forms, which gives the additional
relations
$$dz \cdot dz \,=\,0,\quad dz^* \cdot dz^*\,=\,0,\quad dz^*\cdot
dz \,=\,-q^2dz \cdot dz^*.\eqno(11.3)$$

 Of course, the relations (11.1) -- (11.3) are well known to the specialists
(see, for instance, \cite{M} and the references therein).

\bigskip
\section{Differential calculi: step one}

 We follow the notation of sections 3, 5, 6.

 Consider the linear functionals $\lambda_ \pm \in {\frak h}^*$ given by
$$\lambda_ \pm(H_i) = \pm a_{ij_0},$$
together with the associated generalized Verma modules $V_ \pm(\lambda_
\pm)$. Just as in Section 8, define the "involutions" $*:V_ \pm(\lambda_
\pm) \to V_ \mp(\lambda_ \mp)$ by (8.1) and $*:v_ \pm(\lambda_ \pm) \mapsto
v_ \mp(\lambda_ \mp)$.

 It follows from the definitions that the maps
$$v_+(\lambda_+) \mapsto E_{j_0}v_+(0),\quad v_+(\lambda_+)^*\mapsto
(E_{j_0}v_+(0))^*$$
admit the unique extensions to $U_q{\frak g}$-module morphisms
$$\delta_+:\,V_+(\lambda_+) \to V_+(0), \quad \delta_-:\,V_-(\lambda_-) \to
V_-(0).$$

 Consider the dual graded $U_q{\frak g}$-modules:
$$\bigwedge \nolimits^1({\frak g}_{-1})_q \:=\:\scriptstyle \bigoplus
\limits_{\scriptstyle j \in {\Bbb Z}_+}\displaystyle V_-(\lambda_-)_{-j}^*;
\qquad \bigwedge \nolimits^1(\overline{\frak g}_{-1})_q \:=\:\scriptstyle
\bigoplus \limits_{\scriptstyle j \in {\Bbb Z}_+}\textstyle
V_+(\lambda_+)_j^*.$$
Our definition of the graded components implies that
$$\delta_+V_+(\lambda_+)_j \subset V_+(0)_j;\quad \delta_-V_-(\lambda_-)_j
\subset V_-(0)_j.$$
Now the "adjoint" operators
$\partial=\delta_-^*,\;\overline{\partial}=\delta_+^*$ are well defined and
become $U_q{\frak g}$-module morphisms:
$$\partial:\,{\Bbb C}[{\frak g}_{-1}]_q \to \bigwedge \nolimits^1({\frak
g}_{-1})_q;\quad \overline{\partial}:\,{\Bbb C}[\overline{\frak g}_{-1}]_q
\to \bigwedge \nolimits^1(\overline{\frak g}_{-1})_q.$$

 Evidently, the maps
$$v_ \pm(\lambda_ \pm)\mapsto v_ \pm(0)\otimes v_ \pm(\lambda_ \pm); \quad
v_ \pm(\lambda_ \pm)\mapsto v_ \pm(\lambda_ \pm)\otimes v_ \pm(0)$$
admit the unique extension to $U_q{\frak g}$-module morphisms
$$\Delta_ \pm^L:\,V_ \pm(\lambda_ \pm) \to V_ \pm(0) \otimes V_
\pm(\lambda_\pm), \quad \Delta_ \pm^R:\,V_ \pm(\lambda_ \pm) \to V_
\pm(\lambda_\pm) \otimes V_ \pm(0).$$

 Pass again to the "adjoint" linear operators and observe that they are well
defined and equip $\bigwedge^1({\frak g}_{-1})_q$ with a structure of a
covariant bimodule over ${\Bbb C}[{\frak g}_{-1}]_q$, and
$\bigwedge^1(\overline{\frak g}_{-1})_q$ with a structure of a covariant
bimodule over ${\Bbb C}[\overline{\frak g}_{-1}]_q$. (The covariance here
means that the actions $(\Delta_ \pm^L)^*, (\Delta_ \pm^R)^*$ of ${\Bbb
C}[{\frak g}_{-1}]_q$ and ${\Bbb C}[\overline{\frak g}_{-1}]_q$ respectively
are $U_q{\frak g}$-module morphisms).

\medskip

 Remark. With $\omega \in \bigwedge^1({\frak g}_{-1})_q$ or $\omega \in
\bigwedge^1(\overline{\frak g}_{-1})_q$ one has $1 \cdot \omega\,=\,\omega
\cdot 1\,=\,\omega$, since
$$(\varepsilon \otimes {\rm id})\Delta_ \pm^L(v)\,=\,v, \quad ({\rm
id}\otimes \varepsilon)\Delta_ \pm^R(v)\,=\,v,\quad v \in V_
\pm(\lambda_ \pm).$$

\medskip

 It is easy to show that $\partial$ and $\overline{\partial}$ are
differentiations of the corresponding covariant bimodules:
$$\partial(f_1f_2)\:=\:\partial f_1 \cdot f_2\,+\,f_1 \partial f_2;\quad
f_1,f_2 \in {\Bbb C}[{\frak g}_{-1}]_q,$$
$$\overline{\partial}(f_1f_2)\:=\:\overline{\partial}f_1 \cdot f_2\,+\,f_1
\overline{\partial}f_2;\quad f_1,f_2 \in {\Bbb C}[\overline{\frak
g}_{-1}]_q.$$
For example, to prove the latter inequality, it suffices to
pass in each its part to the adjoint operators
$$V_+(\lambda_+)\,\to\,V_+(0)\otimes V_+(0)$$
and then to apply both operators to the generator $v_+(\lambda_+)$ of the
$U_q{\frak g}$-module $V_+(\lambda_+)$. In both cases one obtains
$$E_{j_0}v_+(0)\otimes v_+(0)\,+\,v_+(0)\otimes E_{j_0}v_+(0).$$

 In conclusion, let us prove one of the equalities (11.1). Another one can
be derived in a similar way.

 It follows from $z^*dz^*\in(V_+(\lambda_+)_2)^*,\; dz^*\cdot z^*\in
(V_+(\lambda_+)_2)^*,\;{\rm dim}\,V_+(\lambda_+)_2\,=\,1$ that $z^*\cdot
dz^*\,=\,{\rm const}\cdot dz^*\cdot z^*$. Thus, it remains to compute the
constant.

 When applying the duality argument, we replace $f(v)$ by $\langle f,v
\rangle$. Compare $\langle z^*dz^*,Ev_+(\lambda_+)\rangle$ and $\langle
dz^*\cdot z^*,Ev_+(\lambda_+)\rangle$.

 Firstly, one has

$\langle z^*dz^*,Ev_+(\lambda_+)\rangle\:=$

$\langle z^*\otimes dz^*,(1 \otimes E\,+\,E \otimes K)(v_+(0)\otimes
v_+(\lambda_+))\rangle\:=$

$\langle z^*\otimes dz^*,(E \otimes K)(v_+(0)\otimes
v_+(\lambda_+))\rangle\:=$

$\langle z^*,Ev_+(0)\rangle \langle dz^*,Kv_+(\lambda_+)\rangle\:=$

$q^2 \langle z^*,Ev_+(0)\rangle \langle dz^*,v_+(\lambda_+)\rangle\:=\:q^2
\langle z^*,Ev_+(0)\rangle^2,$

and secondly

$\langle dz^*\cdot z^*,Ev_+(\lambda_+)\rangle \:=$

$\langle dz^*\otimes z^*,(1 \otimes E\,+\,E \otimes K)(v_+(\lambda_+)\otimes
v_+(0))\rangle \:=$

$\langle dz^*,v_+(\lambda_+)\rangle \langle z^*,Ev_+(0)\rangle \:=\:\langle
z^*,Ev_+(0)\rangle^2.$

Since $\langle z^*,Ev_+(0)\rangle \ne 0$, we obtain finally
$$z^*dz^*\,=\,q^2dz^*\cdot z^*.$$

\bigskip
\section{Differential calculi: step two}

 Consider the $U_q{\frak g}$-module
$$\Omega^{(1,0)}({\frak g}_{-1})_q \,\stackrel{def}{=}\, \bigwedge
\nolimits^1({\frak g}_{-1})_q \otimes {\Bbb C}[\overline{\frak g}_{-1}]_q.$$

 Use the universal R-matrix in the same way as in Section 7 to equip
$\Omega^{(1,0)}({\frak g}_{-1})_q$ with a structure of a covariant bimodule
over ${\rm Pol}({\frak g}_{-1})_q$.

 There is a unique extension of the differentiation $\partial:\,{\Bbb
C}[{\frak g}_{-1}]_q \to \bigwedge^1({\frak g}_{-1})_q$ to a differentiation
$\partial:\,{\rm Pol}({\frak g}_{-1})_q \to \Omega^{(1,0)}({\frak
g}_{-1})_q$ such that $\partial{\Bbb C}[\overline{\frak g}_{-1}]_q \,=\,0$.
Clearly $\partial(f_+ \otimes f_-)\,=\,\partial f_+ \otimes f_-,\quad f_+
\in{\Bbb C}[{\frak g}_{-1}]_q,\;f_- \in{\Bbb C}[\overline{\frak g}_{-1}]_q$,
and $\partial$ is a $U_q{\frak g}$-module morphism.

 Turn to the example ${\frak g}=U_q \frak{sl}_2$. Differentiation of both
sides in (9.10) (with the properties $\partial:1 \mapsto 0,\;\partial:z
\mapsto dz$ being taken into account) yields $z^*\cdot dz-q^2dz \cdot
z^*=0$. This is just one of the relations (11.2).

 Now consider the $U_q{\frak g}$-module $\Omega^{(0,1)}({\frak
g}_{-1})_q \,\stackrel{def}{=}\,{\Bbb C}[{\frak g}_{-1}]_q \otimes
\bigwedge^1(\overline{\frak g}_{-1})_q$ together with the morphism of
$U_q{\frak g}$-modules
$$\overline{\partial}:\,{\rm Pol}({\frak g}_{-1})_q \to
\Omega^{(0,1)}({\frak g}_{-1})_q;\quad \overline{\partial}:\,f_+ \otimes f_-
\mapsto f_+ \otimes \overline{\partial}f_-,$$
where $f_+ \in{\Bbb C}[{\frak g}_{-1}]_q,\;f_- \in{\Bbb C}[\overline{\frak
g}_{-1}]_q$. Just as it was done before, one can equip $\Omega^{(0,1)}({\frak
g}_{-1})_q$ with a structure of a covariant bimodule over ${\rm Pol}({\frak
g}_{-1})_q$ and prove that $\overline{\partial}$ is a differentiation. An
application of $\overline{\partial}$ to both sides of (9.10) gives the
second one of the relations (11.2).

  Finally, set
$$\Omega^1({\frak g}_{-1})_q \,=\,\Omega^{(1,0)}({\frak g}_{-1})_q \,\oplus
\,\Omega^{(0,1)}({\frak g}_{-1})_q,\quad d=\partial+\overline{\partial}.$$

\bigskip
\section{Differential calculi: step three}

 Let $A$ be a unital algebra over ${\Bbb C}$.

 \underline{Definition.} Let $M$ be a bimodule over $A$ and $d:\,A \to M$ a
linear operator. The pair $(M,d)$ is called a differential calculus of order
one if\\
i) $d(a^\prime \cdot a^{\prime \prime})\:=\:da^\prime \cdot a^{\prime
\prime}\,+\,a^\prime \cdot da^{\prime \prime}$,\\
ii) $A \cdot(dA)\cdot A \,=\,M$.

 In the case when $A$ is a covariant algebra, $M$ a covariant bimodule and
$d:\,A \to M$ a $U_q{\frak g}$-module morphism with conditions i), ii) being
satisfied, the pair $(M,d)$ is called a covariant differential calculus of
order one.

 The results expounded in the appendix of this work imply that the five
covariant differential calculi of order one:
$$(\bigwedge \nolimits^1({\frak g}_{-1})_q,\partial),\quad
(\bigwedge\nolimits^1(\overline{\frak g}_{-1})_q,\overline{\partial})$$
$$(\Omega^{(1,0)}({\frak g}_{-1})_q,\partial),\quad (\Omega^{(0,1)}({\frak
g}_{-1})_q,\overline{\partial}),\quad(\Omega^1({\frak g}_{-1})_q,d).$$

 In the sequel we apply to each of those the "algorithm of constructing the
full differential calculus" described in \cite{M}.

 \underline{Definition.} Let $\Omega \,=\,\bigoplus \limits_{n \in {\Bbb
Z}_+}\Omega_n$ be a ${\Bbb Z}_+$-graded algebra and $d$ a linear operator in
$\Omega$ of order one. The pair $(\Omega,d)$ is called a differential graded
algebra if\\
i) $d^2=0$,\\
ii) $d(a^\prime \cdot a^{\prime \prime})\:=\:da^\prime \cdot a^{\prime
\prime}\,+\,(-1)^na^\prime \cdot da^{\prime \prime}, \;a^\prime \in
\Omega_n,\:a^{\prime \prime}\in \Omega$.

 If $\Omega$ is a covariant algebra and $d$ a $U_q{\frak g}$-module
morphism, then under the conditions i) and ii) we call the pair $(\Omega,d)$
a covariant differential graded algebra.

 Let us describe the "algorithm" of construction of the pair $(\Omega,d)$
given the pair $(M,d)$. Let $M_1=dA \subset M$.

 Equip the tensor algebra $T=T(A,M_1)$ with a grading in which ${\rm deg}\,a
\,=\,0,\;{\rm deg}\,m \,=\,1,\;a \in A,\;m \in M_1$. One has $T_0
\:=\:T(A)\:=\:{\Bbb C}\oplus A \oplus A^{\otimes 2}\oplus
\ldots,\;T_{j+1}\,=\,T(A)\otimes M_1 \otimes T_j$.

 There exists a unique operator $d:\,T \to T$ such that\\
i) $d(t_1t_2)\:=\:dt_1 \cdot t_2 \,+\,(-1)^nt_1 \cdot dt_2,\;t_1 \in
T_n,\,t_2 \in T$,\\
ii) $d|_A$ coincides with the differentiation in the initial calculus of
order one,\\
iii) $d|_{M_1}=0$.

 In fact, on $T_0$ we have $d1=0,\;d(a_1 \otimes a_2 \otimes \ldots \otimes
a_k)\,=\,\sum \limits_ja_1 \otimes \ldots \otimes a_{j-1}\otimes da_j
\otimes \ldots \otimes a_k$. From now on we proceed by induction:\\
$d(a \otimes m \otimes t)\:=\:da \otimes m \otimes t\,-\,a \otimes m \otimes
dt,\quad a \in T_0,\;m \in M_1,\;t \in T_j$.

 (Note that $d$ is well defined because of the multilinearity of the
right-hand sides of the above identities in the "indeterminates"
$(a_1,\ldots a_k)$ and $(a,m,t)$ respectively.)

 Consider the least $d$-invariant bilateral ideal $J$ of $T$ which contains
all the elements of the form\\
i) $a_1 \otimes a_2 \,-\,a_1a_2,\;a_1,a_2 \in A$\\
ii) $1 \otimes m-m,\;m \otimes 1-m, \quad m \in M_1$\\
iii) $\{\sum \limits_{ij}a_i^\prime \otimes m_{ij} \otimes a_j^{\prime
\prime}|\;a_i^\prime,a_j^{\prime \prime}\in A,\:m_{ij}\in M_1,\:\sum
\limits_{ij}a_i^\prime m_{ij}a_j^{\prime \prime}\;=\;0\}$\\
(Note that the left hand side of the latter equality is a sum of elements
of the $A$-bimodule $M$.)

 It follows from our construction that $J$ is a graded ideal: $J
\,=\,\bigoplus \limits_j (J \bigcap T_j)$. Furthermore, $J$ is a $U_q{\frak
g}$-submodule of $T$ (due to the covariance of the algebra $A$, the module
$M$ and the order one calculus $(M,d)$).

 Hence the quotient algebra $\Omega=T/J$ with the differential $d_J:\,t+J
\mapsto dt+J$ is a covariant graded differential algebra. It is easy to
show that $A \simeq \Omega_0,\;M \simeq \Omega_1$, and the initial
differential $d:\,A \to M$ "coincides" with the restriction of $d_J$ onto
$\Omega_0$.

 The five order one differential calculi we have already produced lead to
five covariant graded differential algebras
$$(\bigwedge \nolimits_q({\frak g}_{-1}),\partial),\quad(\bigwedge
\nolimits_q(\overline{\frak
g}_{-1}),\overline{\partial}),\quad(\Omega_q^{(*,0)},\partial),\quad
(\Omega_q^{(0,*)},\overline{\partial}),\quad(\Omega_q,d).$$

 In the example ${\frak g}=\frak{sl}_2$ the relations
$\partial^2(z^2)=\overline{\partial}^2((z^*)^2)=\partial(zdz^*-
q^{-2}dz^*z)=0$ imply (11.3).

\bigskip
\section{Holomorphic bundles and Dolbeault
complexes}

 Just as in Section 5, we choose a functional $\mu \in {\frak h}^*$ such
that $m_j=\mu(H_j)\in {\Bbb Z}_+$ for $j \ne j_0$. The linear functional of
this form $\lambda_-\in{\frak h}^*$ was already considered in Section 12.

 Consider a $U_q{\frak g}$-module $V_-(\mu)$ and the associated "graded
dual" module $\Gamma_ \mu$.

 Recall that we use the comultiplication $\Delta^{op}$ to equip
$V_-(0)\otimes V_-(\mu)$ and $V_-(\mu)\otimes V_-(0)$ with a structure of a
$U_q{\frak g}$-module.  Also, the morphisms $$\Delta_L:\:V_-(\mu)\,\to
\,V_-(0)\otimes V_-(\mu);\quad \Delta_L:\:v_-(\mu)\,\mapsto \,v_-(0)\otimes
v_-(\mu),$$ $$\Delta_R:\:V_-(\mu)\,\to \,V_-(\mu)\otimes V_-(0);\quad
\Delta_R:\:v_-(\mu)\,\mapsto \,v_-(\mu)\otimes v_-(0),$$
(together with the adjoint linear maps $\Delta_L^*,\Delta_R^*$) are used to
equip $\Gamma_ \mu$ with a structure of a covariant bimodule over ${\Bbb
C}[{\frak g}_{-1}]_q$:
$$\Delta_L^*:\:{\Bbb C}[{\frak g}_{-1}]_q \otimes \Gamma_ \mu \,\to \,\Gamma_
\mu;\quad \Delta_R^*:\:\Gamma_ \mu \otimes{\Bbb C}[{\frak g}_{-1}]_q \, \to
\,\Gamma_ \mu.$$

 It follows from the properties of the universal R-matrix over
$U_q{\frak g}$ that
$$\sigma R^{-1}v_-(0)\otimes v_-(\mu)\,=\,v_-(\mu)\otimes v_-(0),$$
where $\sigma:\,a \otimes b \mapsto b \otimes a$, and $R^{-1}$ is the
universal R-matrix of the Hopf algebra $U_q{\frak g}^{op}$. Hence $\sigma
R^{-1}\Delta_L=\Delta_R,\quad \Delta_L=R \sigma \Delta_R,$
$$\Delta_L^*=\Delta_R^*\cdot \check{R}.\eqno(15.1)$$
Here $\check{R}:\:{\Bbb C}[{\frak g}_{-1}]_q \otimes \Gamma_ \mu \,\to
\,\Gamma_ \mu \otimes{\Bbb C}[{\frak g}_{-1}]_q,\quad \check{R}=\sigma R$.

 (15.1) shows how to describe the covariant bimodule $\Gamma_ \mu$ in terms
of generators and relations.

 The standard construction (see Section 7) allows one to equip the tensor
product $M_ \mu \,=\,\Gamma_ \mu \otimes{\Bbb C}[\overline{\frak g}_{-1}]_q$
with a structure of a covariant bimodule over ${\rm Pol}({\frak g}_{-1})_q
\,=\,{\Bbb C}[{\frak g}_{-1}]_q \otimes {\Bbb C}[\overline{\frak g}_{-1}]_q$.

 Consider the simplest case ${\frak g}=\frak{sl}_2$. Denote by $\gamma_ \mu$
the lowest weight vector of the $U_q{\frak g}$-module $\Gamma_ \mu$ such
that $\gamma_ \mu(v_-(\mu))=1$. Clearly $m_ \mu=\gamma_ \mu \otimes 1$ is a
generator of a covariant bimodule $M_ \mu$. It is not hard to deduce the
complete list of "commutation" relations using the explicit form of the
universal R-matrix (see Section 9):
$$z \cdot m_ \mu \,=\,q^{-\mu(H)}\cdot m_ \mu \cdot z,\quad z^*\cdot m_
\mu \,=\,q^{\mu(H)}\cdot m_ \mu \cdot z^*.\eqno(15.2)$$

 It is easy to prove that
$$K^{\pm 1}m_ \mu=q^{\pm \mu(H)}m_ \mu,\quad Fm_ \mu=0,\quad Em_
\mu=-q^{1/2} \cdot \frac{1-q^{2 \mu(H)}}{1-q^2}\cdot zm_ \mu.$$
(The last equality follows from the covariance of the bimodule $M_ \mu$
and the relations $Em_ \mu={\rm const}\cdot zm_ \mu,\quad FEm_
\mu=-(EF-FE)m_ \mu=-\frac{\textstyle q^{\mu(H)}-q^{-\mu(H)}}{\textstyle
q-q^{-1}}m_ \mu$.)

 The elements of $M_ \mu$ could be treated as $q$-analogues of smooth
sections of a holomorphic vector bundle. We are interested in differential
forms whose coefficients are such "sections".

 Consider the covariant bimodule $\Omega_{\mu,q}^{(0,*)}\,=\,M_ \mu
\bigotimes_{{\rm Pol}({\frak g}_{-1})_q} \Omega_q^{(0,*)}$ over ${\rm
Pol}({\frak g}_{-1})_q$. Evidently
$$\Omega_{\mu,q}^{(0,*)}\,=\,\Gamma_ \mu \textstyle \bigotimes
\nolimits_{{\Bbb C}[{\frak g}_{-1}]_q} \Omega_q^{(0,*)}\,=\,\Gamma_ \mu
\textstyle \bigotimes \nolimits_{\Bbb C} \bigwedge(\overline{\frak
g}_{-1})_q.\eqno(15.3)$$

 Apply the $U_q{\frak g}$-module morphism $\check{R}:\:{\Bbb C}[{\frak
g}_{-1}]_q \otimes \Gamma_ \mu \,\to \,\Gamma_ \mu \otimes{\Bbb C}[{\frak
g}_{-1}]_q,\quad \check{R}=\sigma R$ derived from the universal R-matrix to
equip $\Omega_{\mu,q}^{(0,*)}$ with a structure of a covariant bimodule
over $\Omega_q^{(0,*)}$. In the example ${\frak g}=\frak{sl}_2$ one can
readily describe this module: the relation list (15.2) should be completed
with one more relation
$$dz^* \cdot m_ \mu \,=\,q^{\mu(H)}\cdot m_ \mu \cdot dz^*.$$

 It follows from (15.3) and $\overline{\partial}{\Bbb C}[{\frak g}_{-1}]_q=0$
that the operator
$$\overline{\partial}_ \mu \:=\:{\rm id}\textstyle \bigotimes
\nolimits_{{\Bbb C}[{\frak g}_{-1}]_q}\overline{\partial}:\,\Gamma_ \mu
\textstyle \bigotimes_{{\Bbb C}[{\frak g}_{-1}]_q}\Omega_q^{(0,*)}\,\to
\,\Gamma_ \mu \textstyle \bigotimes_{{\Bbb C}[{\frak
g}_{-1}]_q}\Omega_q^{(0,*)}.$$

 Certainly, $\Omega_{\mu,q}^{(0,*)}$ is a graded bimodule
$\Omega_{\mu,q}^{(0,*)}\,=\,\bigotimes \limits_j \Omega_{\mu,q}^{(0,j)}$
over $\Omega_q^{(0,*)}$, and $\overline{\partial}_ \mu$ is its
differentiation of order one:
$$\overline{\partial}_ \mu(am)\:=\:(\overline{\partial}a)m \,+\,(-1)^{{\rm
deg}\,a}\cdot a \cdot \overline{\partial}_ \mu m,$$
$$\overline{\partial}_ \mu(ma)\:=\:(\overline{\partial}_ \mu m)a
\,+\,(-1)^{{\rm deg}\,m}\cdot m \cdot \overline{\partial}a$$
for all homogeneous elements $a \in \Omega_q^{(0,*)},\;m \in
\Omega_{\mu,q}^{(0,*)}$.

 Evidently, the differentiation $\overline{\partial}_ \mu$ is determined
unambiguously by its values on generators.

 In the example ${\frak g}=\frak{sl}_2$ for the generator $m_ \mu \in M_ \mu
\hookrightarrow \Omega_{\mu,q}^{(0,*)}$ as above we have:
$$\overline{\partial}m_ \mu=0.$$

 Now pass to the homogeneous components
$$\Omega_{\mu,q}^{(0,j)}\,=\,M_ \mu \textstyle \bigotimes \nolimits_{{\rm
Pol}({\frak g}_{-1})_q} \Omega_q^{(0,j)}\,=\,\Gamma_ \mu \textstyle
\bigotimes \nolimits_{{\Bbb C}[{\frak g}_{-1}]_q}\Omega_q^{(0,j)}$$ of the
graded bimodule $\Omega_{\mu,q}^{(0,*)}$ to obtain the \underline{Dolbeault
complex}
$$0 \to M_ \mu \stackrel{\overline{\partial}_
\mu}{\to}\Omega_{\mu,q}^{(0,1)}\stackrel{\overline{\partial}_
\mu}{\to}\Omega_{\mu,q}^{(0,2)}\stackrel{\overline{\partial}_
\mu}{\to}\ldots.$$
Its terms are the covariant bimodules over ${\rm Pol}({\frak g}_{-1})_q$,
and the differentials are the $U_q{\frak g}$-module morphisms which commute
with the left and the right actions of ${\Bbb C}[{\frak g}_{-1}]_q$.

\bigskip
\section{Conclusion notes.}

 Let us now digress from involutions and differentiations and sketch our
approach to the construction of $q$-analogues of Hermitian symmetric spaces
of non-compact type (one can find more details in sections 2 -- 10).

 Let $q=1$. Evidently, for all $\xi \in {\frak g}_{\pm 1}$ the series
$\exp(\xi)v_ \pm(0)$ converge in some "completed" spaces $\overline{V}_
\pm(0)=\begin{array}[t]{c}\displaystyle \times \\^j \end{array}V_ \pm(0)_j$.
This allows one to elaborate the Harish-Chandra method to produce embeddings
$I_ \pm:\,X \hookrightarrow \overline{V_ \pm(0)}$ of an irreducible
Hermitian symmetric space $X$. The canonical embeddings can be obtained from
$I_ \pm$ by composing them from the right with the projections $\pi_
\pm:\,\overline{V_ \pm(0)}\to V_ \pm(0)_{\pm 1}\simeq {\frak g}_{\pm 1}$.

 Our basic observation is that the topological $U{\frak g}$-modules
$\overline{V_ \pm(0)}$ and hence the subalgebras ${\frak g}_{\pm 1}$ have the
proper quantum analogues \footnote{Remind \cite{CP} that, unlike $U{\frak
g}$, the algebra ${\frak g}$ itself has no "good" quantum analogues.}:
$({\frak g}_{\pm 1})_q \stackrel{def}{=}V_ \pm(0)_ \pm$.

 This allows one to imitate the above Harish-Chandra embeddings $i_ \pm=\pi_
\pm I_ \pm$ for $q \ne 1$.

 There is a different exposition of our construction for $q$-analogues of
bounded symmetric domains and prehomogeneous vector spaces. It provides more
clear interplay between our construction with the approach of V. G. Drinfeld
\cite{Dr1} to quantum groups as well as with the interpretation of the
quantum Weyl group described by S. Z. Levendorski\v{i} and Ya.~S.~Soibelman
\cite{LS}.

 An alternate approach to introducing ${\rm Pol}({\frak g}_{-1})_q$ is in
producing a covariant involutive coalgebra and further passage to the dual
covariant involutive algebra. This approach requires more detailed
exposition of the "duality theory" for $U_q{\frak g}$-module algebras and
$U_q{\frak g}^{op}$-module coalgebras. Specifically, we need to equip our
algebras with the strongest locally compact topologies. The dual coalgebras
are the completions of coalgebras considered in this work above, with
respect to W$^*$-weak topologies, and their tensor products are replaced by
the completed tensor products $\widehat{\otimes}$ (see \cite{KN}). We
describe here the topological covariant $*$-coalgebra dual to ${\rm
Pol}({\frak g}_{-1})_q$. Remind (see section 6) that we replace $U_q{\frak
g}$ by $U_q{\frak g}^{op}$ in tensor products of generalized Verma modules.

 It is easy to show that the vector $v_0=v_-(0)\otimes v_+(0)$ is a
generator of the topological $U_q{\frak g}$-module $\overline{V}_0 \,=\,
\overline{V_-(0)} \widehat{\otimes}\overline{V_+(0)}$. The structure of a
covariant  coalgebra in $\overline{V}_0$ is imposed by introducing a
$U_q{\frak g}$-module morphism
$\overline{\Delta}:\,\overline{V}_0 \to \overline{V}_0
\widehat{\otimes}\overline{V}_0$ given by an application of  a universal
R-matrix:
$$\overline{\Delta}v_0=Rv_0 \otimes v_0.$$
The coassociativity of $\overline{\Delta}$ follows from the
quasitriangularity of the Hopf algebra $U_q{\frak g}$. Impose an involution
in $\overline{V}_0$ by
$$(\xi v_0)^*=(S^{-1}(\xi))^*v_0,\quad \xi \in U_q{\frak g},$$
which already implies
$$(\xi v)^*=(S^{-1}(\xi))^*v^*,\quad \xi \in U_q{\frak g},\; v \in
\overline{V}_0.$$
(Note that (7.1) provides $*$ to be an antilinear coalgebra antihomomorphism
of $\overline{V}_0$).

 Consider the maps
$$\varepsilon_-\otimes {\rm id}:\overline{V}_0 \to \overline{V_+(0)}; \quad
{\rm id}\otimes \varepsilon_+:\overline{V}_0 \to \overline{V_-(0)},$$
with $\varepsilon_ \pm$ being the counits of the coalgebras $\overline{V_
\pm(0)}$.

 It follows from the relations
$$(\varepsilon \otimes{\rm id})(R)=({\rm id}\otimes \varepsilon)(R)=1$$
that these maps are the morphisms of covariant coalgebras (dual to the
embeddings ${\Bbb C}[{\frak g}_{-1}]_q \hookrightarrow {\rm Pol}({\frak
g}_{-1})_q,\quad {\Bbb C}[\overline{\frak g}_{-1}]_q \hookrightarrow {\rm
Pol}({\frak g}_{-1})_q$).

 The relation
$$R(v_0 \otimes v_0)=v_-(0) \otimes R \sigma(v_-(0) \otimes v_+(0)) \otimes
v_+(0)$$
with $\sigma:a \otimes b \mapsto b \otimes a$, demonstrates that the
comultiplication $\overline{\Delta}$ agrees with the multiplication in ${\rm
Pol}({\frak g}_{-1})_q$ introduced in Section 7.

 Finally, let us note that the commutation relations between the elements of
$({\Bbb C}[{\frak g}_{-1}]_q)_{+1}$ and $({\Bbb
C}[\overline{\frak g}_{-1}]_q)_{-1}$ are of degree at most two, as it follows
from the properties of the universal $R$-matrix.

\bigskip
\section*{Appendix. Images of the differentials
$\partial,\overline{\partial}$}

 Consider the coalgebra $V_-(0)$ and the comodule $V_-(\lambda_-)$ and
disregard for a moment their $U_q{\frak g}$-module structures.

\begin{lemma}\hspace{-.5em}. The left comodule $V_-(\lambda_-)$ over the
coalgebra $V_-(0)$ is isomorphic to a direct sum of $m$ copies of $V_-(0)$,
with $m={\rm dim}V_-(\lambda_-)_{-1}$.
\end{lemma}

 {\bf Proof.} Remind the decomposition $w_0=w_0^\prime \cdot u$, with $w_0$
and $w_0^\prime$ being the maximum length elements in the Weyl group of Lie
algebras ${\frak g}$ and ${\frak g}^\prime$ respectively (see Section 5). It
was our agreement to consider only those reduced decompositions of $w_0$
which are given by concatenation of reduced decompositions for $w_0^\prime$
and $u$.

 Choose such a decomposition. Just as in \cite[Proposition 1.7(c)]{CK},
associate to it the bases in $U_q{\frak g}^\prime,\;U_q{\frak g}$ and the
bases of the vector spaces
$$V_-(0)=U_q{\frak g}  v_-(0);\quad V_-(\lambda_-)_{-1}=U_q{\frak
g}^\prime  v_-(\lambda_-);\quad V_-(\lambda_-)=U_q{\frak g}
v_-(\lambda_-).$$
(One can verify that $(V_-(\lambda_-))_{-1}$ is a simple $U_q{\frak
g}^\prime$-module). The above bases are of the form
$$\{\xi_iv_-(0)\},\quad \{\eta_jv_-(\lambda_-)\},\quad \{\xi_i
\eta_jv_-(\lambda_-)\},$$
with $i \in {\Bbb Z}_+,\;j \in \{1,\ldots,m \}$, and $\xi_i \in U_q{\frak
g},\;\eta_j \in U_q{\frak g}^\prime$ can be derived from the bases of
$U_q{\frak g}$ and $U_q{\frak g}^\prime$ described explicitly in \cite{CK}.
It remains to prove that for each $j$ the map
$$\pi_j:\,\xi_iv_-(0)\mapsto \xi_i \eta_jv_-(\lambda_-),\quad i \in {\Bbb
Z}_+$$
is a morphism of left comodules. So
$$\Delta_-^L
\pi_j(\xi_iv_-(0))=\Delta(\xi_i)\Delta_-^L(\eta_jv_-(\lambda_-))=
\Delta(\xi_i)(v_-(0) \otimes \eta_jv_-(\lambda_-))={\rm id}\otimes \pi_j
\Delta(\xi_iv_-(0)).$$

\begin{lemma}\hspace{-.5em}. Let $v \in V_-(0)$. Then $\Delta_-(v)\in
v_-(0)\otimes V_-(0)$ if and only if $v \in {\Bbb C}v_-(0)$.
\end{lemma}

 {\bf Proof.} If $v={\rm const}\cdot v_-(0)$, then
$\Delta_-(v)=v_-(0)\otimes{\rm const}\cdot v_-(0) \in v_-(0)\otimes V_-(0)$.
Conversely, if $\Delta_-(v)=v_-(0)\otimes v_1,\;v_1 \in V_-(0)$, then
$$v=({\rm id}\otimes \varepsilon_-)\Delta_-(v)=({\rm id}\otimes
\varepsilon_-)(v_-(0)\otimes v_1)=\varepsilon_-(v_1)\cdot v_-(0) \in {\Bbb
C}\cdot v_-(0).$$

\begin{lemma}\hspace{-.5em}. Let $v \in V_-(\lambda_-)$. Then
$\Delta_-^L(v)\ \in v_-(0)\otimes V_-(\lambda_-)$ if and only if $v \in
V_-(\lambda_-)_{-1}$.
\end{lemma}

 {\bf Proof.} Let $L=\{v \in V_-(\lambda_-\}|\,\Delta_-^L(v)\in
v_-(0)\otimes V_-(\lambda_-)\}$. Evidently, $L \supset V_-(\lambda_-)_{-1}$,
and by lemmas 1,2, ${\rm dim}L={\rm dim}V_-(\lambda_-)_{-1}$. It follows
that $L=V_-(\lambda_-)_{-1}$.

\begin{remark}\hspace{-.5em}.\end{remark} If $\Delta_-^L(v)=v_-(0)\otimes
v_1$, then $v_1=v$ since $v_1=(\varepsilon \otimes{\rm id})(v_-(0)\otimes
v_1)=(\varepsilon \otimes{\rm id})\Delta_-^L(v)$.

\begin{lemma}\hspace{-.5em}. The restriction of $\delta_-$ onto
$(V_-(\lambda_-))_{-1}$ is an injective linear operator.
\end{lemma}

 {\bf Proof.} This operator is non-zero and is a $U_q{\frak g}$-module
morphism $(V_-(\lambda_-))_{-1}\to V_-(0)$, where $(V_-(\lambda_-))_{-1}$ is
a simple $U_q{\frak g}^\prime$-module.

\begin{proposition}\hspace{-.5em}. ${\Bbb C}[{\frak g}_{-1}]_q \cdot
\partial{\Bbb C}[{\frak g}_{-1}]_q=\bigwedge^1({\frak g}_{-1})_q$.
\end{proposition}

 {\bf Proof.} Assume the contrary. Let $V^\prime=\{v \in
V_-(\lambda_-)|\,\langle f_1 \partial f_2,v \rangle=0,\,\forall f_1,f_2
\in{\Bbb C}[{\frak g}_{-1}]_q \}$. Then $V^\prime=\scriptstyle \bigoplus
\limits_{\scriptstyle i \in {\Bbb Z}_+}\textstyle V_{-i}^\prime \ne 0$. It
follows from the definitions that $\Delta_-^L(V^\prime)\subset V_-(0)\otimes
V^\prime$.

 Let $i^\prime$ be the least such $i \in {\Bbb Z}_+$ that $V_{-i}^\prime
\ne 0$. We have
$$\Delta_-^L(V_{-i^\prime}^\prime)\subset({\Bbb C}v_-(0)\scriptstyle
\bigoplus \scriptstyle \bigoplus \limits_{\scriptstyle k>0}\textstyle
V_-(0)_{-k})\otimes V^\prime \subset v_-(0)\otimes
V_-(\lambda_-)+(\scriptstyle \bigoplus \limits_{\scriptstyle k>0} \textstyle
V_-(0)_{-k})\otimes V^\prime.$$

 On the other hand $((\scriptstyle \bigoplus \limits_{\scriptstyle
k>0}\textstyle V_-(0)_{-k})\otimes V^\prime)_{-i^\prime}=0$, and hence
$\Delta_-^L(V_{-i^\prime}^\prime)\subset v_-(0)\otimes V_-(\lambda_-)$. Now
we deduce from lemma 3 that $V_{-i^\prime}^\prime=V_-(\lambda_-)_{-1}\bigcap
V^\prime$.

 Let $v^\prime \in V_-(\lambda_-)_{-1}\bigcap V^\prime$. It follows from the
definition of $V^\prime$ that $({\rm id}\otimes
\delta_-)\Delta^L_-(v^\prime)=0$. On the other hand, by lemma 3 and remark 4
we observe that $\Delta_L(v^\prime)=v(0)\otimes v^\prime$. Hence $({\rm
id}\otimes \delta_-)(v(0)\otimes v^\prime)=0$. That is,
$\delta_-(v^\prime)=0$ and hence, by lemma 5, $v^\prime=0$.

 Thus we have proved that $V_{-i^\prime}^\prime=0$ which makes a
contradiction to the contrary of proposition 6.

\begin{remark}\hspace{-.5em}.\end{remark} One can prove in a similar way
that
$$\overline{\partial}{\Bbb C}[\overline{\frak g}_{-1}]_q \cdot{\Bbb
C}[\overline{\frak g}_{-1}]_q={\Bbb C}[\overline{\frak g}_{-1}]_q
\cdot \overline{\partial}{\Bbb C}[\overline{\frak g}_{-1}]_q=\bigwedge
\nolimits^1(\overline{\frak g}_{-1})_q,$$
$$\partial{\Bbb C}[{\frak g}_{-1}]_q \cdot{\Bbb C}[{\frak
g}_{-1}]_q=\bigwedge \nolimits^1({\frak g}_{-1})_q.$$

\bigskip

\begin{center}

\end{center}

\begin{thebibliography}{99}

\bibitem{A} E. Abe. Hopf Algebras, Cambridge Univ. Press, Cambridge, 1980.

\bibitem{BR} P. N. Bopp and H. Rubenthaler, {\it Fonction z\^eta associ\'ee
\`a la s\'erie principale sph\'erique de certain espaces symm\'etriques},
Ann. Sci. \'Ec. Norm. Sup., $4^e$ serie, t. 26, 1993, p. 701 - 745.

\bibitem{CP} V. Chari and A. Pressley. A Guide to Quantum Groups, Cambridge
Univ. Press, 1995.

\bibitem{CK} C. de Concini, V. Kac, {\it Representations of quantum groups
at roots of 1}, in Operator Algebras, Unitary Representations, Enveloping
Algebras and Invariant Theory, A. Connes, M. Duflo, A. Joseph, R.
Rentschler (eds.), 1990, Birkhauser, Boston, p. 471 - 506.

\bibitem{Di} J. Dixmier. Les $C^*$-alg\`ebres et leur Repr\'esentations.
Paris, Gauthier-Villars, 1964.

\bibitem{Dr1} V. G. Drinfeld, {\it Quantum groups}, in Proceedings of the
International Congress of Mathematicians, Berkeley, 1986, A. M. Gleason
(ed.), 1987, American Mathematical Society, Providence, R. I., 798 - 820.

\bibitem{Dr2} V. G. Drinfeld, {\it On almost commutative Hopf algebras},
Leningrad Math. J., {\bf 1}, (1990), 321 - 432.

\bibitem{He} S. Helgason. Differential Geometry and Symmetric Spaces, Acad.
Press, N.-Y. -- London, 1962.

\bibitem{Hu} J. E. Humphreys. Reflection Groups and Coxeter groups.
Cambridge Univ. Press, 1990.

\bibitem{JS} A. Joyal and R. Street, {\it Braided tensor categories}, Adv.
in Math., {\bf 102} (1993), 20 - 78.

\bibitem{K} C. Kassel. Quantum Groups. Springer-Verlag,
N.Y.--Berlin--Heidelberg, 1995.

\bibitem{KL} S. Klimek and A. Lesniewski, {\it A two-parameter quantum
deformation of the unit disc}, J. Funct. Anal. {\bf 115}, (1993), 1 - 23.

\bibitem{KN} J. L. Kelley, I. Namioka. Linear Topological Spaces, Van
Nostrand Inc., Princeton N.Y. -- London, 1963.

\bibitem{KRR} S. Khoroshkin, A. Radul, V. Rubtsov, {\it A family of Poisson
structures on compact Hermitian symmetric spaces}, Comm. Math. Phys. {\bf
152}, (1993), 299 - 316.

\bibitem{L} G. Lustig, {\it Quantum groups at roots of 1}, Geometriae
Dedicata {\bf 35}, (1990), 89 - 114.

\bibitem{LS} S. Z. Levendorski\v{i} and Ya. S. Soibelman, {\it Some
applications of the quantum Weyl group}, J. Geom. Phys. {\bf 7}, (1990),
241 - 254.

\bibitem{M} G. Maltsiniotis, {\it Le langage des espaces et des groupes
quantiques}, Commun. Math. Phys., {\bf 151}, (1993), 275 - 302.

\bibitem{NN} G. Nagy, A. Nica. {\it On the "quantum disc" and a
"non-commutative circle"}, in: Algebraic Methods on Operator Theory, R. E.
Curto, P. E. T. J\/orgensen (eds.), Birkhauser, Boston, 1994, p. 276 - 290.

\bibitem{R} H. Rubenthaler, {\it Les paires duales dans les alg\`ebres de
Lie r\'eductives}, Asterisque, v. 219, 1994.

\bibitem{S} J. P. Serre. Complex Semisimple Algebras.
Berlin--Heidelberg--New York, Springer, 1987.

\bibitem{SiV} S. Sinel'shchikov and L. Vaksman, {\it Hidden symmetry of the
differential calculus on the quantum matrix space}, to appear in J. Phys. A.

\bibitem{SoV} Ya. S. Soibelman, L. L. Vaksman, {\it On some problems in the
theory of quantum groups}, in Representation Theory and Dynamical Systems,
A. M. Vershik (ed.), Advances in Soviet Mathematics, {\bf 9}, (1990),
American Mathematical Society, Providence, R. I., 3 - 55.

\bibitem{VSh} S. Sinel'shchikov, D. Shklyarov, and L. Vaksman. On function
theory in the quantum disc: integral representations. Preprint, 1997, q-alg.

\bibitem{VSo} L. L. Vaksman and Ya. S. Soibelman, {\it Algebra of functions
on the quantum group $SU(2)$}, Funct. Anal. Appl., {\bf 22} (1988), No \ 3,
170 - 181.

\bibitem{W} S. L. Woronowicz, {\it Compact matrix pseudogroups}, Comm. Math.
Phys., {\bf 111} (1987), 613 - 665.

\end{thebibliography}
\end{document}